\documentclass[10pt]{amsart}
\usepackage{amsmath, amsthm, amsfonts, amssymb}
\usepackage{ifthen, enumerate}

\newtheorem{thm}{Theorem}[section]
\newtheorem{cor}[thm]{Corollary}
\newtheorem{prop}[thm]{Proposition}

\newtheorem{rem}[thm]{Remark}
\newtheorem{rems}[thm]{Remarks}
\theoremstyle{definition}

\newtheorem{example}[thm]{Example}
\newtheorem{lemma}[thm]{Lemma}
 \numberwithin{equation}{section}
 \newtheorem{assum}[thm]{Assumptions}
 \def\cya{{\mathcal C}_{{\mathcal P},\alpha}}
 \def\cyb{{\mathcal C}_{{\mathcal P},\beta}}
 \def\bya{{\mathcal B}_{{\mathcal P},\alpha}}
 \def\byb{{\mathcal B}_{{\mathcal P},\beta}}
 \def\ea{\mathfrak a}
 \def\eb{\mathfrak b}
 \def\ef{\mathfrak f}
 \def\eg{\mathfrak g}
 \def\eh{\mathfrak h}
 
 \def\eu{\mathfrak u}
 \def\mv{\mathsf{v}}
 \def\me{\mathsf{e}}

\title[coupled parabolic systems]{Qualitative properties of coupled parabolic systems of evolution equations}
\author{Stefano Cardanobile}
\author{Delio Mugnolo}
\address{Institut f\"ur Angewandte Analysis, Universit\"at Ulm, Helmholtzstra{\ss}e 18, D-89081 Ulm, Germany}
\email{stefano.cardanobile@uni-ulm.de, delio.mugnolo@uni-ulm.de}


\thanks{The authors wish to express their gratitude to Wolfgang Arendt and Markus Biegert (Ulm) for helpful discussions.
Part of this paper has been written while the second author was visiting Abdelaziz Rhandi at the University Cadi Ayyad of Marrakech (Morocco), whom he warmly thanks for his hospitality.}

\begin{document}

\maketitle

\begin{abstract}
We apply functional analytical and variational methods in order to study well-posedness and qualitative properties of evolution equations on product Hilbert spaces. To this aim we introduce an algebraic formalism for matrices of sesquilinear mappings. We apply our results to parabolic problems of different nature: a coupled diffusive system arising in neurobiology, a strongly damped wave equation, a heat equation with dynamic boundary conditions.
\end{abstract}

\section{Introduction}

Our aim in this paper is to extend some ideas and techniques introduced by R. Nagel in~\cite{Na89} to investigate systems of linear partial differential equations by means of operator matrices. In his paper, the basic intuition was that a linear algebraic formalism also for matrices of unbounded operators may help to discuss well-posedness and spectral issues in analogy to standard matrix analysis. Instead of dealing with general operator matrices, we introduce suitable matrices of sesquilinear mappings and then investigate well-posedness of differential systems by the elegant theory of sesquilinear forms on Hilbert spaces. In order to fix the ideas we first present our setting.

\begin{assum}\label{assum1}
Throughout this paper we impose the following, for $i,j=1,\ldots,m$.
\begin{enumerate}[(i)]
\item $H_i,V_i$ are complex Hilbert spaces such that $V_i$ is continuously and densely embedded in $H_i$.
\item $a_{ij}:V_j\times V_i\to \mathbb C$ are sesquilinear mappings, i.e., mappings that are linear in the first and antilinear in the second variable.
\end{enumerate}
\end{assum}

We always denote by ${\mathcal H}:=\prod_{i=1}^m H_i$ and ${\mathcal V}:=\prod_{i=1}^m V_i$ the product Hilbert spaces endowed with the canonical scalar products
$$(\ef\mid \eg)_{\mathcal H}:=\sum_{i=1}^m (f_i\mid g_i)_{H_i}\qquad\hbox{and}\qquad (\ef\mid \eg)_{\mathcal V}:=\sum_{i=1}^m (f_i\mid g_i)_{V_i}$$
for $ \ef,\eg\in{\mathcal H}$ and $\ef,\eg\in{\mathcal V}$, respectively. Here and in the following we write $\ef$ for $(f_1,\ldots,f_m)^\top$, and likewise for $\eg,\eh$, etc.
Observe that $\mathcal V$ is continuously and densely embedded in $\mathcal H$.

We introduce a densely defined, sesquilinear form $\ea$ on $\mathcal V$ defined by
\begin{equation}\label{formmatr}{\ea}(\ef,\eg):=\sum_{i,j=1}^m a_{ij}(f_j,g_i),\qquad \ef,\eg\in {\mathcal V}.\end{equation}
Since $\overline{\mathcal V}={\mathcal H}$, there exists a canonical operator $\mathcal A$ associated with $\ea$ given by
\begin{eqnarray*}
D({\mathcal A})&:=&\{\ef \in {\mathcal V} : \exists \eg\in{\mathcal H} \hbox{ s.t. } \ea(\ef,\eh)=(\eg\mid \eh)_{\mathcal H}\hbox{ for all }\eh\in {\mathcal V}\},\\
{\mathcal A} \ef&:=&-\eg.
\end{eqnarray*}
Similarly, we can associate with each mapping $a_{ij}:V_j\times V_i\to \mathbb C$ an operator $A_{ij}$ from $H_j$ to $H_i$ by
\begin{eqnarray*}
D(A_{ij})&:=&\{f_j \in V_j : \exists g_i\in H_i \hbox{ s.t. } a_{ij}(f_j,h_i)=(g_i\mid h_i)_{H_i}\hbox{ for all }h_i\in V_i\},\\
A_{ij}f_j&:=&-g_i.
\end{eqnarray*}

Let us now briefly discuss the special case where for $i\not=j$ the mappings $a_{ij}$ can be extended continuously to the whole product space $H_j\times H_i$, so that each operator $A_{ij}$ is bounded from $H_j$ to $H_i$. Then it is possible to identify the operator $\mathcal A$ associated with $\ea$ with some ease. 

\begin{prop}\label{identif}
Assume that for $i\not=j$ the sesquilinear mappings $a_{ij}$ are continuous on $H_j\times H_i$.
Then the operator $\mathcal A$ associated with $\ea$ has diagonal domain $D({\mathcal A}):=\prod_{i=1}^m D(A_{ii})$ and it is given by
$$\label{domin2}{\mathcal A}\ef:=\left(\sum_{i=1}^m A_{1i}f_i,\ldots,\sum_{i=1}^m A_{mi}f_i\right)^\top,\qquad \ef\in D({\mathcal A}),$$
or rather, in matrix form
\begin{equation}\label{formale}
{\mathcal A}=\begin{pmatrix} A_{11} & \cdots & A_{1m}\\
\vdots & \ddots & \vdots\\
A_{m1} & \cdots & A_{mm}
 \end{pmatrix}.
\end{equation}
\end{prop}

\begin{proof}
Let $\ef\in\mathcal V$ be such that there exists a vector $\eg\in\mathcal H$ satisfying $\ea(\ef,\eh)=(\eg\mid \eh)_{\mathcal H}$ for all $\eh\in\mathcal V$. Observe that if the form $a_{ij}$ is associated with $A_{ij}\in{\mathcal L}(H_j, H_i)$, then 
$$\ea (\ef,\eh)=\sum_{i,j=1}^m a_{ij}(f_j,h_i)=\sum_{i=1}^m a_{ii}(f_i, h_i)-\sum_{i\not= j} (A_{ij}f_j\mid h_i)_{H_i}.$$
On the other hand,
$$\ea (\ef,\eh)=(\eg\mid \eh)_{\mathcal H}=\sum_{i=1}^m (g_i\mid h_i)_{H_i}.$$
In particular, considering vectors of the form $\eh=(0,\ldots,h,\ldots,0)^\top\in\mathcal H$, we see that 
$a_{ii}(f_i,h)-\sum_{j\not= i} (A_{ij}f_j\mid g)_{H_i}=(g_i\mid h)_{H_i}$, i.e.,
$$a_{ii}(f_i,h)=\Big(g_i+\sum_{j\not= i} A_{ij}f_j\mid h\Big)_{H_i}=:(\tilde{g}_i\mid h_i)_{H_i}$$ 
holds for some $\tilde{g}_i\in H_i$ and all $h\in H_i$. It follows from the definition of the operator associated with $a_{ii}$ that $f\in D(A_{ii})$, and $A_{ii}f_i=-\tilde{g}_i=-g_i-\sum_{j\not= i} A_{ij}f_j$. Summing up, $\sum_{j=1}^m A_{ij}f_j =-g_i$ for all $i=1\ldots,m$. It can be proven likewise that the converse inclusion holds.
\end{proof}

Such a casual interpretation of an entrywise interplay between form $\ea$ and operator matrix $\mathcal A$ is not always justified. In Section~\ref{epha} we consider the case of a form whose associated operator is of the type described above although the assumptions of Theorem~\ref{identif} are not satisfied. However, it may as well be that $D({\mathcal A})$ is not a product space, or furthermore that some $A_{ij}=0$ although $a_{ij}\not\equiv 0$, cf. Sections~\ref{csop}--\ref{dbc}, respectively.
Still, we keep the above identification as a heuristic motivation for characterizing generator properties of $\mathcal A$, as well as some features of the generated semigroups, by means of the individual mappings $a_{ij}$.  In a certain sense, this is the same target pursued in~\cite{Na89}. In the spirit of Nagel's article, in most of our results we deduce properties of $\ea$ from individual conditions on $a_{ij}$.

We believe that there are good reasons to develop a matrix theory for forms. First, we show in Section~2 that whole classes of differential problems fit our framework, including evolution equations that do not look like systems of parabolic equations. Furthermore, our matrix formalism allows us to check simple, linear algebraic properties of finite-dimensional matrices, instead of dealing with complicated infinite-dimensional problems. 

Another reason to treat systems by means of sesquilinear forms is that invariance of subsets of the state space can be obtained by a criterion due to E.M. Ouhabaz, cf.~\cite[Thm.~2.2]{Ou04}. We extensively use it in order to investigate invariance properties of sets that, in our opinion, are particularly relevant for systems of coupled evolution equations.

Finally, we emphasize that the setting in~\cite{Na89} is more general than ours. In fact, Nagel considers the case of $C_0$-semigroups, whereas T. Kato has shown that only analytic, quasi-contractive semigroups (and not even all of them) can be generated by operator associated with forms. On the other hand, in~\cite{Na89} only very mild forms of coupling could be treated, cf. the results in~\cite[\S~3]{Na89}: in particular, no well-posedness result was proved for the case where all off-diagonal operators $A_{ij}$ in~\eqref{formale} are ``as unbounded" as the diagonal ones. We consider some possibile applications in Section~4. A further class of systems that fit our theory is given by coupled diffusion--ODE problems of FitzHugh--Nagumo type, see e.g.~\cite{CM07}.

Our results should be compared with those obtained by H. Amann in~\cite{Am84} and E.M. Ouhabaz in~\cite{Ou99} for parabolic problems with state space $L^p(\Omega,H)$, where $H$ is an arbitrary Hilbert space. Well-posedness for a general class of coupled diffusion systems has been discussed in~\cite{AN93}. Finally, let us mention that a rich and elegant theory for operator matrices (both with diagonal and non-diagonal domain), in particular concerning asymptotics of semigroups, has been developed by K.-J. Engel in~\cite{En08}.

\section{Matrices of forms}\label{sectionwellp}

For given Hilbert spaces $V,H$ such that $V\hookrightarrow H$ and numbers $M,\omega\geq 0$ and $\alpha>0$, a sesquilinear form $a:V\times V\to \mathbb C$ is said to be \emph{continuous with constant $M$} and \emph{$H$-elliptic with constants $(\alpha,\omega)$} if
$$|a(u,v)|\leq M \| u\|_{V} \|v\|_{V},\qquad u,v\in V,$$
and
$${\rm Re}a(u,u) \geq \alpha \| u\|^2_V-\omega \| u\|^2_H,\qquad u\in V,$$
respectively. It is said to be \emph{coercive with constant $\alpha$} if it is $H$-elliptic with constants $(\alpha,0)$, and \emph{accretive} if ${\rm Re}a(u,u)\geq 0$ for all $u\in V$. 

By Kato's form characterization of sectorial operators, cf.~\cite[\S~5.3.4]{Ar04}, the operator $A$ associated with $a$ generates an analytic semigroup $(e^{z{A}})_{z\in \Sigma_\theta}$ of angle $\theta\in(0,\frac{\pi}{2}]$ such that $\| e^{z {A}}\|_{{\mathcal L}({\mathcal H})}\leq e^{\omega |z|}$, $z\in\Sigma_{\theta}$, for some $\omega\in\mathbb R$, if and only if $a$ is densely defined, continuous, and $H$-elliptic; such a semigroup is contractive if and only if $a$ is accretive. Thus, we are interested in continuity and ellipticity properties for the form $\ea$ introduced in Section~1.

To begin with, we recall the following perturbation lemma, cf.~\cite[Lemma~2.1]{Mu07b}. For $\alpha\in [0,1)$ we denote by $H_\alpha$ any interpolation space between $V$ and $H$, i.e., any linear space $V \hookrightarrow H_\alpha \hookrightarrow H$ that verifies the interpolation inequality
\begin{equation*}
\|f\|_{H_\alpha}\leq M_\alpha\|f\|_V^\alpha \|f\|_H^{1-\alpha},\qquad f\in V.
\end{equation*}

\begin{lemma}\label{perturbd}
Let $a:V\times V\to\mathbb C$ be a sesquilinear mapping. Let $\alpha\in [0,1)$ such that $a_1:V\times H_\alpha\to \mathbb C$ and $a_2:H_\alpha\times V\to\mathbb C$ are continuous sesquilinear mappings. Then $a$ is $H$-elliptic if and only if $a+a_1+a_2:V\times V\to\mathbb C$ is $H$-elliptic.
\end{lemma}

Observe that the optimal $H$-ellipticity constants of $a$ and $a+a_1+a_2$ is in general different.

In the following immediate consequence of Lemma~\ref{perturbd}, $H_{i\alpha}$ denotes an interpolation space between $V_i$ and $H_i$.

\begin{cor}\label{perturb0}
Let $a_{ii}$ be $H_i$-elliptic for all $i=1,\ldots,m$. Let all off-diagonal sesquilinear mappings $a_{ij}$ be continuous on $V_j\times H_{i\alpha}$ or on $H_{j\alpha}\times V_{i}$. Then also the form matrix $\ea$ is ${\mathcal H}$-elliptic.
\end{cor}

\begin{prop}\label{continpure}
The form $\ea$ is continuous if and only if 
\begin{itemize}
\item for $i=1,\ldots,m$ the forms $a_{ii}$ are continuous with constant $M_{ii}$, and
\item for $i,j=1,\ldots,m$, $i\not=j$, the forms $a_{ij}$ are continuous in the following sense: there exist $\alpha_{ij}\leq 0$ and $\omega_{ij} \in\mathbb R$ such that
\begin{equation}\label{contaij}
|a_{ij}(f,g)|\leq -\alpha_{ij} \| f\|_{V_j} \|g\|_{V_i} + \omega_{ij}\| f\|_{H_j} \|g\|_{H_i},\qquad \hbox{for all }(f,g)\in V_j\times V_i.
\end{equation}
\end{itemize}
In this case the continuity estimates 
$$
|\ea(\ef,\eg)| \leq (\|\mathbf M\|+\|{\mathbf \Omega_0\|}{\mathfrak e^2})\|\ef\|_\mathcal V\|\eg\|_\mathcal V.
$$
holds, where the scalar matrices $\bf M$ and $\bf\Omega_0$ are given by
$${\mathbf M}:=\begin{pmatrix}
M_{11} & -\alpha_{12} & \cdots & -\alpha_{1m}\\
-\alpha_{21} & M_{22} & & -\alpha_{2m}\\
\vdots &&\ddots&\vdots\\
-\alpha_{m1} & -\alpha_{m2} & \cdots & M_{mm}
 \end{pmatrix},\quad 
{\mathbf \Omega_0}:=\begin{pmatrix}
0 & |\omega_{12}| & \cdots & |\omega_{1m}|\\
|\omega_{21}| &0& & |\omega_{2m}|\\
\vdots &&\ddots&\vdots\\
|\omega_{m1}| & |\omega_{m2}| & \cdots &0
 \end{pmatrix}.$$
\end{prop}

Here and in the following, $\mathfrak e$ stands for the norm of the canonical injection of $\mathcal V$ into $\mathcal H$.

\begin{proof}
Let $\ef,\eg\in{\mathcal V}$ and observe that by assumption
\begin{eqnarray*}
|\ea(\ef,\eg)| &\leq& \sum_{i=1}^m |a_{ii}(f_j,g_i)| +\sum_{i\not=j} |a_{ij}(f_j,g_i)| \\
&\leq&\left(\sum_{i=1}^m M_{ii} \|f_i\|_{V_i }\|g_i\|_{V_i}-\sum_{i\not=j}^m \alpha_{ij} \|f_j\|_{V_j}\|g_i\|_{V_i}\right)+\sum_{i\not=j}^m |\omega_{ij}| \|f_j\|_{ H_j}\|g_i\|_{H_i} \\
&\leq &\|{\bf M}\| \| \ef\|_{\mathcal V} \|\eg\|_{\mathcal V}+\|{\bf \Omega}_0\|\; \|\ef\|_\mathcal H \|\eg\|_\mathcal H,
\end{eqnarray*}
where the last step follows from the Cauchy--Schwarz inequality.
This shows one implication. Assume now that $\ea$ is continuous but $a_{i_0j_0}$ is not, for some $i_0,j_0$. Consider sequences $(u_k)_{k\in\mathbb N} \subset V_{j_0}$ and $(v_k)_{k\in\mathbb N}\subset V_{i_0}$ such that $\|u_k\|_{V_{j_0}}=\|v_{k}\|_{V_{i_0}}=1$ for all $k \in \mathbb N$, but $\lim_{k\to\infty}|a_{i_0 j_0}(u_k,v_k)| = \infty$. Define 
$$
{\mathfrak u}_k:=\begin{pmatrix}
0\\\vdots\\ u_k\\ \vdots\\ 0\end{pmatrix}\leftarrow j_0^{\rm th} \hbox{ row},\qquad
{\mathfrak v}_k:=\begin{pmatrix}
0\\\vdots\\ v_k\\ \vdots\\ 0\end{pmatrix}\leftarrow i_0^{\rm th} \hbox{ row},\qquad k\in{\mathbb N},
$$ 
One sees that $\|\mathfrak u_k\|_\mathcal V=\|\mathfrak v_k\|_\mathcal V=1$ for all $k\in\mathbb N$, and there holds $|\ea(\mathfrak u_k, \mathfrak v_k)|=|a_{i_0j_0}(u_k,v_k)|$, a contradiction to the continuity of $\ea$.
\end{proof}

In the following we focus on the case where off-diagonal mappings $a_{ij}$ are actually unbounded on $H_j\times H_i$, since Corollary~\ref{perturb0} and Proposition~\ref{continpure} already allow us to discuss parabolic problems whose associated forms have off-diagonal bounded entries with respect to some interpolation space.

We recall that a scalar $m\times m$ matrix $M=(m_{ij})$ is called \emph{positive} (resp. \emph{negative}) \emph{semidefinite} if there exists $\mu\geq 0$ such that $(M\xi\cdot \xi)\geq \mu |\xi|^2$ (resp., $(M\xi\cdot \xi)\leq -\mu|\xi|^2$)
for all $\xi\in{\mathbb C}^m$. Further, $M$ is called \emph{positive definite} (resp. \emph{negative definite}) if it is positive (resp. negative) semidefinite  and $\mu$ can be chosen $>0$.

\begin{prop}\label{basic}
The following assertions hold for the densely defined form $\ea$.
\begin{enumerate}[(1)]
\item If the form $\mathfrak a$ is $\mathcal H$-elliptic with constants $(\alpha,\omega)$, then for all $i=1,\ldots,m$ the forms $a_{ii}$ are ${H_i}$-elliptic with constants $(\alpha,\omega)$, too.
\item Conversely, assume $a_{ii}$ to be $H_i$-elliptic with constants $(\alpha_{ii},\omega_{ii})$, $i=1,\ldots,m$. Let~\eqref{contaij} hold, and assume the matrix ${\bf A}:=(\alpha_{ij})_{1\leq i,j\leq m}$ to be positive definite with constant $\alpha>0$. Then the form $\ea$ is $\mathcal H$-elliptic with constants $(\alpha,\|{\bf\Omega}\|)$, where ${\bf \Omega}=(\omega_{ij})_{1\leq i,j\leq m}$.
\item If $\ea$ is accretive, then all $a_{ii}$ are accretive, $i=1,\ldots,m$.
\item Let~\eqref{contaij} hold and $a_{ii}$ be accretive, $i=1,\ldots,m$. If the matrices ${\mathbf A}_0:={\mathbf A}-{\rm diag}({\mathbf A})$ and
${\mathbf \Omega}_0:={\mathbf  \Omega}-{\rm diag}({\mathbf  \Omega})$
are positive and negative semidefinite, respectively, then $\ea$ is accretive.
\end{enumerate}
\end{prop}

\begin{proof}
Assertions (1) and (3) can be checked in a way that is similar to that used in the proof of Proposition~\ref{continpure}, by considering vectors of the form $\ef=(0,\ldots,f,\ldots,0)^\top$.	

In order to prove (2) and (4), let now~\eqref{contaij} hold. Then for all $\ef\in{\mathcal V}$
\begin{eqnarray*}
{\rm Re} \ea(\ef,\ef)&=&{\rm Re} \sum_{i=1}^m a_{ii}(f_i,f_i)+ {\rm Re} \sum_{j \not= i}^m a_{ij}(f_j,f_i)\\
&\geq& \sum_{i,j=1}^m \alpha_{ij}\|f_j\|_{V_j}\|f_i\|_{V_i} - \sum_{i,j=1}^m \omega_{ij}\|f_j\|_{H_j}\|f_i\|_{H_i}\\
%
&\geq& \alpha \|\ef\|^2_\mathcal V - \|{\bf \Omega}\|\;	\|\ef\|^2_\mathcal H.
\end{eqnarray*}
Likewise in (4), if all forms $a_{ii}$ are accretive, then 
\begin{eqnarray*}
{\rm Re} \ea(\ef,\ef)&=&{\rm Re} \sum_{i=1}^m a_{ii}(f_i,f_i)+ {\rm Re} \sum_{j\not= i}^m a_{ij}(f_j,f_i)\geq {\rm Re} \sum_{j\not= i}^m a_{ij}(f_j,f_i)\\
&\geq& \sum_{j\not= i}^m \alpha_{ij}\| f_j\|_{V_j} \| f_i\|_{V_i} -\sum_{j\not= i}^m \omega_{ij}\| f_j\|_{H_j} \| f_i\|_{H_i}.
\end{eqnarray*}
This shows that ${\rm Re} \ea(\ef,\ef)\geq 0$ if ${\bf A}_0$ and ${\bf \Omega}_0$ are positive and negative semidefinite, respectively.
\end{proof}

\begin{rem}\label{gersh}
Assume $a_{ii}$ to be $H_i$-elliptic (resp., coercive), $i=1,\ldots,m$. If furthermore
\begin{equation}\label{gersapplic}
\alpha_{ii}> \sum _{k\not= i} \frac{|\alpha_{ik}+\alpha_{ki}|}{2}, \qquad k=1,\ldots,m,
\end{equation}
then it follows from Gershgorin's circle theorem that $\sigma(\frac{{\bf A}+{\bf A^*}}{2})$ is contained in the open right half plain of $\mathbb C$. Since the coercivity of $\bf A$ is equivalent to the strict positive definiteness of $\frac{{\bf A}+{\bf A^*}}{2}$, it follows that~\eqref{gersapplic} is a sufficient condition for $\ea$ to be $\mathcal H-$elliptic. In particular, we can \emph{always} obtain well-posedness by suitably weakening the strength of the internal coupling of the system, i.e., by letting individual parameters $\alpha_{ij}\to 0$.
If in particular $\omega_{ij}=0$, $i,j=1,\ldots,m$, then Gershgorin's circle theorem also yield a threshold beyond which the semigroup associated with $\ea$ is exponentially stable. 
\end{rem}

The following is motivated by Proposition~\ref{basic}.



\begin{assum}\label{assum2}
In the remainder of the paper we impose the following, for $i,j=1,\ldots,m$, $j\not=i$	.
\begin{enumerate}[(i)]
\item $a_{ii}$ is continuous with constant $M_{i}$ and $H_i$-elliptic with constants $(\alpha_i,\omega_i)$.
\item $a_{ij}$ satisfies~\eqref{contaij} for constants $\omega_{ij},\alpha_{ij}$  such that ${\bf A}=(\alpha_{ij})_{1\leq i,j\leq m}$ is positive definite.
\end{enumerate}
\end{assum}

By~\cite[Prop.~1.51 and Thm.~1.52]{Ou04} we obtain well-posedness for the abstract Cauchy problem
\begin{equation}\tag{ACP}
\left\{\begin{array}{rcll}
\dot{\eu}(t)&=&\mathcal A\mathcal\eu(t),\qquad &t\geq 0,\\
\eu(0)&=&\eu_0.
 \end{array}\right.
\end{equation} 

\begin{thm}\label{semig}
The operator $\mathcal A$ associated with $\ea$ generates on $\mathcal H$ an analytic semigroup of angle $\frac{\pi}{2}-\arctan (\| {\bf M}\| +\|{\bf \Omega}_0\| {\mathfrak e}^2)$. This semigroup is compact if and only if $V_i$ is compactly embedded in $H_i$ for all $i=1,\ldots,m$. It is uniformly exponentially stable if for $\omega_{ij}=0$, $i,j=1\ldots,m$, $\bf A$ is positive definite.
\end{thm}

The estimate on the analyticity angle obtained in Theorem~\ref{semig} can often be improved.

\begin{prop}\label{cosine}
The following assertions hold.
\begin{enumerate}[(1)]
\item Assume that there exists $M\geq 0$ such that for all $f\in V_j$ and $g\in V_i$ one has
\begin{enumerate}[(i)]
\item $|{\rm Im}a_{ii} (f,f)|\leq M \| f\|_{V_i} \| f\|_{H_i}$ for all $i=1,\ldots,m,$ and moreover 
\item \begin{itemize}
\item either $|{\rm Im}(a_{ij}(f,g)+{a_{ji}(g,f)})|\leq M\| f\|_{V_j}\| g\|_{H_i}$, for all $i,j=1,\ldots,m \hbox{ s.t. } i<j,$ 
\item or $|{\rm Im}(a_{ij}(f,g)+{a_{ji}(g,f)})|\leq M\| g\|_{H_j}\| f\|_{V_i}$, {for all} $i,j=1,\ldots,m \hbox{ s.t. } i<j.$
\end{itemize}
\end{enumerate}
Then the operator $\mathcal A$ associated with $\ea$ generates a cosine operator function with associated phase space ${\mathcal V}\times {\mathcal H}$. In particular, $\mathcal A$ generates an analytic semigroup of angle $\frac{\pi}{2}$ on $\mathcal H$.
	
\item Conversely, if $\mathcal A$ generates a cosine operator function, then for all $i=1,\ldots,m$ also the operator $A_{ii}$ associated with $a_{ii}$ generates a cosine operator function.
\end{enumerate}
\end{prop}

\begin{proof}
Under the assumptions in (1), we have
\begin{eqnarray*}
| {\rm Im}\ea (\ef,\ef)|&\leq& \left|\sum_{i=1}^m {\rm Im}\Big(a_{ii}(f_i,f_i)\right|+\left|\sum_{i\not=j} {\rm Im}\Big(a_{ij}(f_j,f_i)\Big)\right|\\
&\leq &\left|\sum_{i=1}^m {\rm Im}\Big(a_{ii}(f_i,f_i)\right|+\left|\sum_{i<j} {\rm Im}\Big(a_{ij}(f_j,f_i)+a_{ji}(f_i,f_j)\Big)\right|\\
&\leq &\sum_{i=1}^m |{\rm Im}\Big(a_{ii}(f_i,f_i)|+\sum_{i<j} |{\rm Im}\Big(a_{ij}(f_j,f_i)+a_{ji}(f_i,f_j)\Big)	|\\
&\leq &\left\{\begin{array}{l}
M \| f_i\|_{V_i} \| f_i\|_{H_i} + M\sum_{i<j}^m \| f_j\|_{V_j}\| f_i\|_{H_i}\\
M \| f_i\|_{V_i} \| f_i\|_{H_i} + M\sum_{i<j}^m \| f_i\|_{V_i}\| f_j\|_{H_j}
 \end{array}\right.\\
&\leq& \tilde{M}\| \ef\|_{\mathcal V} \| \ef\|_{\mathcal H}
\end{eqnarray*}
for some constant $\tilde{M}\geq 0$. Applying a result due to Crouzeix--Haase in the version presented in~\cite[p.~204]{Ha06}, one obtains that $\mathcal A$ generates a cosine operator function with associated phase space ${\mathcal V}\times {\mathcal H}$, hence also an analytic semigroup of angle $\frac{\pi}{2}$ on $\mathcal H$ by~\cite[Thm.~3.14.17]{ABHN01}.

By the above mentioned result of Crouzeix--Haase, $\mathcal A$ generates a cosine operator function if and only if there exists an equivalent scalar product $(\!(\cdot |\cdot)\!)$ on $\mathcal H$ such that the numerical range 
$$W(\ea):=\{\ea(\eu,\eu)\in{\mathbb C}: \eu\in {\mathcal V} \hbox{ and } (\!( \eu|\eu)\!)=1\}$$
lies in a parabola, cf.~\cite[\S~5.6.6]{Ar04}. In order to prove (2) consider the equivalent scalar product on $\mathcal H$ with respect to which $W(\ea)$ lies in a parabola. Such a scalar product on $\mathcal H$ induces an equivalent scalar product on $H_i$, too, and therefore also the numerical range $W(a_{ii})$ lies in the same parabola, for all $i=1,\ldots,m$. Thus, $A_{ii}$ generates a cosine operator function by Crouzeix--Haase's result.
\end{proof}


\section{Averaging and invariance properties}\label{sectionsymm}

Having investigated the well-posedness of the Problem (ACP), we turn our attention to qualitative properties of the semigroup associated with the form $\ea$ introduced in \eqref{formmatr}, which can be described by means of the invariance of suitable subsets of the state space. In this Section we still impose Assumptions~\ref{assum1} and \ref{assum2}.

The following result characterizes the invariance of product subspaces. It is a direct consequence of Corollary~\ref{lemmino} and we omit its easy proof.

\begin{prop}\label{productinv}
Let $Y_i$ be closed subspaces of the Hilbert spaces $H_i$ for each $i=1,\ldots,m$ and denote by $P_i$ the corresponding orthogonal projections. Then the subspace $\mathcal Y:=\prod_{i=1}^m Y_i$ is invariant under the action of the semigroup $(e^{t\ea})_{t\geq0}$ if and only if
\begin{itemize}
\item $P_j V_j \subset V_j$ for all $j=1,\ldots,m$, and
\item $a_{ij}(f,g)=0$ for all $f \in Y_j\cap V_j, g\in Y_i^\perp\cap V_i$ and all $i,j=1,\ldots,m$.
\end{itemize}
\end{prop}

We can characterize invariance of a special class of subspaces of ${\mathcal H}=\prod_{i=1}^mH_i$ that cannot be represented as a Cartesian product. In~\cite{CMN07} we have also discussed in detail the interplay between invariance of such kind of subspaces and the notion of symmetries of a physical system.

\begin{thm}\label{linainvar}
Let $(X,\mu)$ a $\sigma$-finite meausre space such that $H_1=\ldots=H_m=L^2(X)$, i.e., $\mathcal H= L^2(X)^m\simeq L^2(X;{\mathbb C^m})$. Assume furthermore that $V_1=\ldots=V_m$ and the form $\ea$ to be accretive.
Consider an orthogonal projection $K=(\kappa_{ij})_{1\leq i,j\leq m}\in M_m({\mathbb C})$ and define the operator $${\mathcal P}\ef:=K\ef=\left(\sum_{j=1}^m \kappa_{1j}f_j,\ldots,\sum_{j=1}^m \kappa_{mj}f_j\right)^\top,\qquad \ef\in{\mathcal H}.$$ 
Let $\mv_i=(v_{i1},\ldots,v_{im})^\top$, $i=1,\ldots,m$, be an orthonormal basis of eigenvectors of ${\mathbb C}^m$ such that $\mv_1,\ldots,\mv_r$ are eigenvectors of $K$ associated with eigenvalue 1, and
$\mv_{r+1},\ldots,\mv_m$ are eigenvectors of $K$ associated with eigenvalue 0.
\begin{enumerate}[(1)]
\item The following assertion are equivalent.
\begin{enumerate}[(a)]
\item The semigroup $(e^{ta})_{t\geq 0}$ leaves invariant the closed subsets $\cya:=\left\{\ef\in {\mathcal H}:\|\ef-{\mathcal P}\ef \|\leq \alpha \right\}$ for some/all $\alpha\geq0$;
\item for all $\ef\in \mathcal V$, $\eg\in{\rm ker}(I-{\mathcal P})\cap {\mathcal V}$, and $\eh\in{\rm ker}({\mathcal P})\cap {\mathcal V}$ there holds ${\mathcal P}\ef\in \mathcal V$ and $\ea(\eg, \eh)= 0$;
\item for all $\eg\in \mathcal V$ there holds
\begin{equation}\label{autovk1}
\sum _{i,j=1}^m\sum _{\ell=1}^r\sum _{k=r+1}^m v_{\ell j}\overline{v_{k i}}a_{ij}(g_\ell,g_k)=0.
\end{equation}
\end{enumerate}
\item Furthermore, the following assertions are also equivalent.
\begin{enumerate}[(a')]
\item The semigroup $(e^{ta})_{t\geq 0}$ leaves invariant the closed subsets $\bya:=\left\{\ef\in {\mathcal H}:\|{\mathcal P}\ef \|\leq \alpha \right\}$ for some/all $\alpha\geq0$;
\item for all $\ef\in \mathcal V$, $\eg\in{\rm ker}(I-{\mathcal P})\cap {\mathcal V}$, and $\eh\in{\rm ker}({\mathcal P})\cap {\mathcal V}$ there holds ${\mathcal P}\ef\in \mathcal V$ and $\ea(\eh, \eg)= 0$;
\item for all $\eg\in\mathcal V$ there holds
\begin{equation}\label{autovk2}
\sum _{i,j=1}^m\sum _{\ell=1}^r\sum _{k=r+1}^m v_{ j\ell}\overline{v_{ ik}}a_{ij}(g_\ell,g_k)=0.
\end{equation}
\end{enumerate}
\end{enumerate}
\end{thm}

\begin{proof}
We only show that (1.a)--(1.c) are equivalent, the proof of the equivalences in (2) being analogous.

First of all, observe that the linear operator $\mathcal P$ is an orthogonal projection on $\mathcal H$: in fact, it is a contraction that satisfies ${\mathcal P}={\mathcal P}^2$, due to the analogous properties of the matrix $K$. The equivalence of (1.a)--(1.b) is then a direct consequence of Corollary~\ref{lemmino}. In order to prove that (1.b) is equivalent to (1.c), observe that each coordinate of $\mathcal P \ef$ is a linear combination of $f_1,\ldots,f_m$, thus again a vector of $V$: thus, $\mathcal P \ef\in\mathcal V$. Consider now the projection $K$ in its Jordan nomal form to see that its eigenvalues are 0 and/or 1, i.e., $\sigma(K) \subset \{0,1\}$, and that it is diagonalizable. Thus, it is always possible to find $\mv_1,\ldots,\mv_m$ with the required properties. 

Let $\ef\in \mathcal H$ and decompose the vector $\ef(x)\in{\mathbb C}^m$ as $\ef(x)=\sum _{j=1}^m \lambda^\ef_i(x) \mv_i$. Observe that $\lambda^\ef_1,\ldots,\lambda^\ef_m\in L^2(X)$, and in fact $\lambda^\ef_1,\ldots,\lambda^\ef_m\in V$ if $\ef \in\mathcal V$. Moreover, \begin{eqnarray*}
\mathcal P\ef(x)&=& \left(\sum_{i,j=1}^m \lambda^\ef_{i}(x) \kappa_{1j} v_{ij},\ldots,\sum_{i,j=1}^m \lambda^\ef_{i}(x)\kappa_{mj} v_{ij}\right)^\top \\
&=&\left(\sum_{i=1}^m \lambda^\ef_{i}(x)\sum_{j=1}^m \kappa_{1j} v_{ij},\ldots,\sum_{i=1}^m \lambda^\ef_{i}(x)\sum_{j=1}^m\kappa_{mj} v_{ij}\right)^\top\\
&=&\left(\sum_{i=1}^m \lambda^\ef_{i}(x)(K\mv_i)_1,\ldots,\sum_{i=1}^m \lambda^\ef_{i}(x)(K\mv_i)_m\right)^\top,
\end{eqnarray*}
holds for $\mu$-almost every $x\in X$. Since now $K\mv_i=\mv_i$, $i=1,\ldots,r$ and $K\mv_i=0$, $i=r+1,\ldots,m$, there holds
$$
\mathcal P\ef= \left(\sum_{i=1}^r \lambda^\ef_{i} v_{i1},\ldots,\sum_{i=1}^r \lambda^\ef_{i} v_{im}\right)^\top, \qquad
(I-\mathcal P)\ef= \left(\sum_{i=r+1}^m \lambda^\ef_{i} v_{i1},\ldots,\sum_{i=r+1}^m \lambda^\ef_{i} v_{im}\right)^\top, \qquad \mu\mbox{-a.e.}
$$
Accordingly, there holds
$$
\ker (I-\mathcal P)=\left\{\eg\in{\mathcal H}: \exists \ef\in\mathcal H \hbox{ s.t. } \eg= \sum _{i=1}^r\lambda^\ef_i \mv_i\right\}=\left\{\eg\in{\mathcal H}: \exists \lambda_1,\ldots,\lambda_r\in H \hbox{ s.t. } \eg= \sum _{i=1}^r\lambda_i \mv_i\right\}
$$
as well as
$$
\ker (\mathcal P)=\left\{\eh\in{\mathcal H}: \exists \ef'\in\mathcal H \hbox{ s.t. } \eh= \sum _{i=r+1}^m\lambda^{\ef'}_i \mv_i\right\}=\left\{\eg\in{\mathcal H}: \exists \lambda_{r+1},\ldots,\lambda_m\in H \hbox{ s.t. } \eg= \sum _{i=r+1}^m\lambda_i \mv_i\right\},
$$
so that
$$
\ker (I-\mathcal P)\cap {\mathcal V}=\left\{\eg\in{\mathcal H}: \exists \lambda_1,\ldots,\lambda_r\in V \hbox{ s.t. } \eg= \sum _{i=1}^r\lambda_i \mv_i\right\}
$$
and
$$
\ker (\mathcal P)\cap {\mathcal V}=\left\{\eg\in{\mathcal H}: \exists \lambda_{r+1},\ldots,\lambda_m\in V \hbox{ s.t. } \eg= \sum _{i=r+1}^m\lambda_i \mv_i\right\}.
$$

We are finally in the position to prove the equivalence of (1.b) and (1.c). In fact, let $\eg\in\mathcal V$ and decompose $\eg=\eg_1\oplus \eg_2$, with $\eg_1\in\ker (I-\mathcal P)\cap {\mathcal V}$ and $\eg_2\in\ker (\mathcal P)\cap {\mathcal V}$. Then by (1.b) one has 
$$
	0=\ea(\eg_1,\eg)=\sum _{i,j=1}^m a_{ij} \big(\sum _{\ell=1}^r \lambda_{\ell} v_{\ell j},\sum _{k=r+1}^m \lambda_k v_{ki}\big)= \sum _{i,j=1}^m\sum _{\ell=1}^r\sum _{k=r+1}^m v_{\ell j}\overline{v_{k i}}a_{ij}(\lambda_\ell,\lambda_k).$$
Similarly if $\eg=\sum_{i=1}^r \lambda_i \mv_i\in{\rm ker}(I-{\mathcal P})\cap {\mathcal V}$ and $\eh=\sum_{i=r+1}^m \lambda_i \mv_i \in{\rm ker}({\mathcal P})\cap {\mathcal V}$, then
$$
\ea(\eg,\eh)=\sum _{i,j=1}^m a_{ij} \big(\sum _{\ell=1}^r \lambda_{\ell} v_{\ell j},\sum _{k=r+1}^m \lambda_k v_{ki}\big)= \sum _{i,j=1}^m\sum _{\ell=1}^r\sum _{k=r+1}^m v_{\ell j}\overline{v_{k i}}a_{ij}(\lambda_\ell,\lambda_k)=0.$$
This concludes the proof.
\end{proof}

Let us consider again the abstract Cauchy problem (ACP) introduced in Section~2. In the case of a system whose state space is $L^2(X)\times L^2(X)$, it seems interesting to consider under which assumptions initial conditions that are ``in phase" (i.e., such that $u_{01}=u_{02}$) give rise to solutions to (ACP) that are in phase as well (i.e., such that $u_1(t)=u_2(t)$), cf. Remark~\ref{blbekbes} below. A natural generalization of this problem is discussed in the following.

\begin{example}\label{invarmed}
Let $(X,\mu)$ be a $\sigma$-finite measure space. Consider a Hilbert space $V$ such that $V\hookrightarrow H:=L^2(X)$ and ${\mathcal H}:=H^m$, ${\mathcal V}:=V^m$. Consider an accretive form $\ea$ and a linear operator ${\mathcal P}$  defined by
$${\mathcal P}\ef:=\left(\sum_{j=1}^m \frac{f_i}{m},\ldots, \sum_{j=1}^m \frac{f_i}{m}\right)^\top,\qquad \ef\in{\mathcal H}.$$ 
Then the semigroup $(e^{ta})_{t\geq 0}$ leaves invariant closed subsets $\cya:=\left\{\ef\in {\mathcal H}:\|\ef-{\mathcal P}\ef \|\leq \alpha \right\}$ for some/all $\alpha\geq0$ if and only if
\begin{equation}\label{primainvarcond}
\sum _{i,j=1}^m a_{ij}(g,h_i)=0, \qquad\hbox{for all } g,h_1,\ldots,h_m\in V \hbox{ s.t. }\sum_{i=1}^m h_i\equiv 0.
\end{equation}
In fact, ${\mathcal P}\ef=K\ef$ for all $\ef\in\mathcal H$, where $K=(\kappa_{ij})_{1\leq i,j\leq m}$ with $\kappa_{ij}=\frac{1}{m}$. One checks that $K$ is an orthogonal projection and $\ef\in{\rm ker}({\mathcal P})$ if and only if $\sum_{i=1}^m f_i(x)=0$ for $\mu$-a.e. $x\in X$, while $\ef\in{\rm ker}(I-{\mathcal P})$ if and only if $f_i(x)=f_j(x)=:f(x)$ for $\mu$-a.e. $x\in X$ and all $i,j=1,\ldots,m$. Thus we deduce by Theorem~\ref{linainvar}.(1)	 that
$(e^{ta})_{t\geq 0}$ leaves invariant closed subsets $\cya:=\left\{\ef\in {\mathcal H}:\|\ef-{\mathcal P}\ef \|\leq \alpha \right\}$ for some/all $\alpha\geq0$ if and only if for all $\eg\in\ker(I-{\mathcal P})\cap \mathcal V$ and all $\eh\in\ker({\mathcal P})\cap \mathcal V$ there holds $\ea(\eg,\eh)=0$, i.e., if and only if
$$\sum_{i,j=1}^m a_{ij}(g_j,h_i)=\sum_{i,j=1}^m a_{ij}(g,h_i)=0\qquad\hbox{for all }g,h_1\ldots,h_m\in V \hbox{ s.t. }\sum_{i=1}^m h_i \equiv 0.$$
Likewise one can see that the semigroup $(e^{ta})_{t\geq 0}$ leaves invariant the closed subsets $\bya:=\left\{\ef\in {\mathcal H}:\|{\mathcal P}\ef \|\leq \alpha \right\}$ for some/all $\alpha\geq 0$ if and only if
\begin{equation}\label{primainvarcond2}
\sum _{i,j=1}^m a_{ij}(g_j,h)=0, \qquad\hbox{for all } g_1,\ldots,g_m,h\in V \hbox{ s.t. }\sum_{i=1}^m g_i\equiv 0.
\end{equation}
In the special case of ${\mathcal H}=L^2(X)\times L^2(X)$, we have thus characterized under which assumptions initial conditions ``in counterphase" give rise to solutions to (ACP) that are in counterphase, too.
\end{example}


Theorem~\ref{linainvar} also allows to study invariance of subsystems.

\begin{example}\label{connec}
Let $(X,\mu)$ be a $\sigma$-finite measure space. Consider a Hilbert space $V$ such that $V\hookrightarrow H:=L^2(X)$ and ${\mathcal H}:=H^m$, ${\mathcal V}:=V^m$. Let the form $\ea$ be accretive and consider the linear operator ${\mathcal P}$ defined by ${\mathcal P}\ef:=(f_1,\ldots,f_{m_0},0,\ldots,0)^\top$, for some $m_0\in\{2,\ldots,m-1\}$. Then the semigroup $(e^{t\ea})_{t\geq 0}$ leaves invariant the closed convex set 
$$\cya:=\prod_{i=1}^{m_0} H_i\times \prod_{i=m_0+1}^m \{ f\in H_i: \| f\|_{H_i}\leq \alpha\}$$
for some/all $\alpha\geq 0$ if and only if the forms $a_{ij}=0$ for all $i=m_0+1,\ldots,m$ and all $j=1,\ldots, m_0$.

Indeed, $\mathcal P$ defined above is the orthogonal projection of $\mathcal H$ onto $\prod_{i=1}^{m_0} H_{i}\times \{0\}^{m-m_0}$. One see that ${\mathcal P}\ef\in \mathcal V$ for all $\ef\in\mathcal V$. In order to apply Theorem~\ref{linainvar} let 
$$
K:=\begin{pmatrix} I&0\\0&0 \end{pmatrix},
$$
where $I$ is the identity $m_0\times m_0$ matrix. Denote with $\me_i$, $i=1,\ldots,m$ the vectors of the canonical basis of $\mathbb C^m$ and observe that $\me_1,\ldots,\me_{m_0}$ are eigenvectors of $K$ associated with eigenvalue 1, whereas $\me_{m_0+1},\ldots,\me_{m}$ are eigenvectors associated with eigenvalue 0. This implies that Theorem~\ref{linainvar}.(1) applies with $r:=m_0$ and $v_{ij}:=\delta_{ij}$, $i,j=1,\ldots,m$ if and only if for all $\eg \in\mathcal V$
$$\sum _{i,j=1}^m\sum _{\ell=1}^{m-1}\sum _{k=m}^m \delta_{\ell j}\overline{\delta_{k i}}a_{ij}(g_\ell,g_k)	
=\sum_{i=m_0+1}^m\sum _{j=1}^{m_0} a_{ij}(g_j,g_i)=0,$$
This condition is satisfied if and only if $a_{ij}=0$ for all $i=m_0+1,\ldots,m$ and all $j=1,\ldots,m_0$.
\end{example}

In the remaining of this section we prove results that can only be formulated whenever our Hilbert state space $\mathcal H$ is an $L^2$-space. Thus, we throughout assume that $H_i=L^2(X_i)$ for a $\sigma$-finite measure space  $(X_i,\mu_i)$, $i=1,\ldots,m$. Accordingly, we can identify ${\mathcal H}$ with $L^2(X)$, where $(X,\mu)$ is a suitable $\sigma$-finite measure space such that $\mu=\mu_1\oplus\ldots\oplus\mu_m$.

\begin{thm}\label{contractive}
Let $H_i=L^2(X_i)$. Then the following assertions hold.
\begin{enumerate}[(1)]
\item The semigroup $(e^{t\ea})_{t\geq 0}$ is real, i.e., it leaves invariant the subset of real-valued functions in $\mathcal H$, if and only if
\begin{itemize}
\item $f\in V_i\Longrightarrow {\rm Re}f\in V_i,\hbox{ and } a_{ii}({\rm Re}f,{\rm Im}f)\in{\mathbb R}\hbox{ for all } i=1,\ldots,m,\hbox{ and}$
\item $a_{ij}(f,g)\in{\mathbb R}\hbox{ for all real-valued }f\in V_j,\; g\in V_i,\; i,j=1,\ldots,m$, $i\not=j$,
\end{itemize}
\item The semigroup $(e^{t\ea})_{t\geq 0}$ is positive, i.e., it leaves invariant the positive cone of $\mathcal H$, if and only if it is real and moreover
\begin{itemize}
\item $f\in V_i\Longrightarrow ({\rm Re})^+f\in V_i$, and $a_{ii}(({\rm Re}f)^+,({\rm Re}f)^-)\leq 0$ for all $ i=1,\ldots,m$, and
\item $a_{ij}(f,g)\leq 0$ for all $0\leq f\in V_j$ and $0\leq g \in V_i$, $i,j=1,\ldots,m$, $i\not=j$.
\end{itemize}
\item Let $(e^{t\ea})_{t\geq 0}$ be positive. Consider another densely defined, continuous, ${\mathcal H}$-elliptic sesquilinear form $\eb:=\sum_{i,j=1}^m b_{ij}:{\mathcal W}\times{\mathcal W}\to {\mathbb C}$, ${\mathcal W}=\prod_{i=1}^m W_i$. Then $(e^{t\ea})_{t\geq 0}$ dominates $(e^{t\eb})_{t\geq 0}$ in the sense of positive semigroups if and only if
\begin{itemize}
\item $W_i$ is an ideal of $V_i$ for all $i=1,\ldots,m$,
\item ${\rm Re} b_{ii}(f,g) \geq a_{ii}(|f|,|g|)$ for all $f,g\in V_i$ such that $f\overline{g}\geq 0,\; i=1,\ldots,m$, and
\item $|{\rm Re} b_{ij}(f,g)| \leq -a_{ij}(|f|,|g|)$ for all $f \in V_j, g\in V_i,\; i=1,\ldots,m.$
\end{itemize}
\end{enumerate}
\end{thm}

\begin{proof}
By Proposition~\ref{basic} the form $\ea $ is densely defined, continuous, and $\mathcal H$-elliptic. Thus, by~\cite[Prop. 2.5 and Thm. 2.6]{Ou04}, and taking into account a rescaling argument, the semigroup $(e^{t \ea })_{t\geq 0}$ is real, positive, and dominating $(e^{t\eb})_{t\geq 0}$, respectively, if and only if
\begin{enumerate}[(i)]
\item ${\ef}\in \mathcal V \Rightarrow {\rm Re}\ef\in \mathcal V \hbox{ and } { \ea }({\rm Re}{\ef},{\rm Im}{\ef})\in\mathbb{R}$, and
\item ${\ef}\in \mathcal V\Rightarrow ({\rm Re} {\ef})^+ \in \mathcal V,\; {\ea}({\rm Re}\ef,{\rm Im}\ef)\in{\mathbb R}$, and $ \ea (({\rm Re})^+,({\rm Re})^-)\leq 0$,
\item ${\mathcal W}$ is an ideal of ${\mathcal V}$ and ${\rm Re} \eb (\ef,\eg)\geq \ea (|\ef|,|\eg|)$ for all $f,g\in \mathcal W$ such that $\ef\overline{\eg}\geq 0$,
\end{enumerate}
respectively. 
First, let (i) hold. If $\ef\in \mathcal V$, then ${\rm Re}\; \ef$ and also ${\rm Re}\; \ef^+\in \mathcal V$. Then, by considering vectors of the form $\ef=(0,\ldots,f,\ldots,0)^\top$ one sees that $a_{i_0 i_0}({\rm Re}{f},{\rm Im}{f})\in\mathbb{R}$ for all $i_0=1,\ldots,m$. Take now $i_0\not= j_0$ and let us show that $a_{i_0 j_0}(f,g)\in{\mathbb R}$ for all real-valued $f\in V_j,\; g\in V_i$. Construct a vector $\ef$ so that all its coordinates besides the $i_0^{\rm th}$ and the $j_0^{\rm th}$ ones vanish, and let its $i_0^{\rm th}$ coordinate agree with $i	 g$ and its $j_0^{\rm th}$ coordinate agree with $f$. Then, it follows that $a_{i_0 j_0}(f,g)=\ea({\rm Re}\ef,{\rm Im} \ef)\in{\mathbb R}$.

Likewise, let (ii) hold. We show that the conditions on $a_{ii}$ and $a_{ij}$ in (2) are satisfied. Again by considering vectors of the form $\ef=(0,\ldots,f,\ldots,0)^\top$ one sees that $a_{i_0 i_0}(({\rm Re}{f})^+,({\rm Re}{f})^-)\leq 0$ for all $i_0=1,\ldots,m$. Take now $i_0\not= j_0$ and let $0\leq f\in V_j,\; 0\leq g\in V_i$. Construct a vector $\ef$ so that all its coordinates besides the $i_0^{\rm th}$ and the $j_0^{\rm th}$ ones vanish, and let its $i_0^{\rm th}$ coordinate agree with $-g$ and its $j_0^{\rm th}$ coordinate agree with $f$. Then, it follows that $a_{i_0 j_0}(f,g)=\ea(({\rm Re}\ef)^+,({\rm Re} \ef)^-)\leq 0$. It is easy to convince oneself that the converse implications in (1) and (2) hold, too.


Finally, let (iii) hold. Since all of its factor spaces $W_i$ are ideal of the spaces $V_i$, $i=1,\ldots,m$, it is clear that ${\mathcal W}$ is an ideal of ${\mathcal V}$. In the following, we denote by
\begin{equation}\label{fjhj}
{\ef}^{j_0}:=\begin{pmatrix}
0\\\vdots\\ f\\ \vdots\\ 0\end{pmatrix}\leftarrow j_0^{\rm th} \hbox{ row}\qquad\hbox{and}\qquad
{\eg}^{i_0}:=\begin{pmatrix}
0\\\vdots\\ g\\ \vdots\\ 0\end{pmatrix}\leftarrow i_0^{\rm th} \hbox{ row},
\end{equation}
vectors in $\mathcal V$, for given $f\in V_{j_0}$ and $g\in V_{i_0}$.

Let now $i_0=j_0$ and $f\overline{g}\geq 0$, so that $\ef^{i_0}\overline{\eg}^{i_0}\geq 0$. Computing ${\rm Re} b_{i_0i_0}(f,g)={\rm Re} \eb (\ef,\eg)\geq \ea(|\ef|,|\eg|)= a_{i_0i_0}(|f|,|g|)$
shows that the second condition holds. For $i_0\not=j_0$, let $f \in V_{j_0}$ and $g\in V_{i_0}$, so that $\ef^{j_0}\overline{\eg^{i_0}}=0=(-\ef^{j_0})\overline{\eg^{i_0}}$. Then, 
$$\pm{\rm Re} b_{i_0j_0}(f,g)={\rm Re} b_{i_0j_0}(\pm f,g)={\rm Re} \eb (\pm \ef^{j_0},\eg^{i_0})\geq \ea(|\ef^{j_0}|,|\eg^{i_0}|)= a_{i_0j_0}(|f|,|g|),$$
thus proving that the third condition is necessary. To check the converse implication let $\ef,\eg \in \mathcal W$ and compute 
${\rm Re} \eb(\ef,\eg) = {\rm Re} \sum _{i,j=1}^m b_{ij}(f_j,g_i)\geq\sum _{i,j=1}^m a_{ij}(|f_j|,|g_i|)=\ea(|\ef|,|\eg|).$
\end{proof}

As a direct consequence of Theorem~\ref{contractive}.(3) we state the following.

\begin{cor}\label{domindiag}
Let the semigroup $(e^{t\ea})_{t\geq0}$ be positive, and assume $a_{ij}(f,g)\leq 0$ for all $0\leq f \in V_j$ and $0\leq g\in V_i$, $i,j=1,\ldots,m$. Let $\ea_0:=\sum _{i=1}^ma_{ii}$. Then $(e^{t\ea})_{t\geq0}$ dominates $(e^{t\ea_0})_{t\geq0}$.
\end{cor}

If $(X,\mu)$ is a $\sigma$-finite measure space and a semigroup $(T_2(t))_{t\geq 0}$ is contractive in both $L^2(X)$ and $L^\infty(X)$, then one sees by standard interpolation results that the semigroup extrapolates to a family $(T_p(t))_{t\geq 0}$ of $C_0$-semigroups in all spaces $L^p(X)$, $p>2$. Such a family is consistent in the sense that $T_p(t)f=T_q(t)f$ for all $f\in L^p(X)\cap L^q(X)$. This motivates the following.

\begin{thm}\label{xinfty}
Let $H_i:=L^2(X_i)$. Assume $\ea$ to be accretive. Then $(e^{t\ea})_{t\geq 0}$ is \emph{$L^\infty$-contractive}, i.e., it leaves invariant the unit ball of $L^\infty(X)$, if and only if for all $i=1,\ldots,m$ there holds
\begin{enumerate}[(i)]
\item $f\in V_i\Longrightarrow(1\wedge |f|){\rm sign}f\in V_i$ and
\item $\sum_{j\not= i} |a_{ij}(f_j,(|f_i|-1)^+{\rm sign}f_i)|\leq {\rm Re}a_{ii}((1\wedge|f_i|){\rm sign}f_i,(|f_i|-1)^+{\rm sign} f_i)$ for all $\ef\in {\mathcal V}\cap {\mathcal C}^\infty_i,$ where the sets ${\mathcal C}^\infty_i$ are defined in \eqref{inftyball}
\end{enumerate}
In particular, all semigroups $(e^{ta_{ii}})_{t\geq 0}$ are $L^\infty$-contractive if so is $(e^{t\ea})_{t\geq 0}$. 
\end{thm}

Here, ${\rm sign}f$ denotes the generalized (complex-valued) sign function defined by
$$({\rm sign} f)(x):=\left\{
\begin{array}{ll}
\frac{f(x)}{\vert f(x)\vert}\quad &\hbox{if } f(x)\not=0,\\
0 &\hbox{if } f(x)=0.
\end{array}
\right.$$
Moreover, we denote by $B^\infty_X$ the unit ball of $L^\infty(X)$ and by ${\mathcal C}^\infty_i$ the set 
\begin{equation}\label{inftyball}
{\mathcal C}^\infty_i:=B^\infty_{X_1}\times \ldots B^\infty_{X_{i-1}}\times L^2(X_i)\times B^\infty_{X_{i+1}}\times\ldots\times B^\infty_{X_m},\qquad i=1,\ldots,m.
\end{equation}

\begin{proof}
By~\cite[Thm.~2.14]{Ou04} the semigroup is $L^\infty$-contractive if and only if 
\begin{equation}\label{ouhcon}
\ef\in \mathcal V \Rightarrow (1\wedge \vert {\ef}\vert){\rm sign} \ef \in \mathcal V\;\hbox{and}\;{\rm Re}{\ea}((1\wedge \vert \ef\vert){\rm sign} \ef,(\vert \ef\vert-1)^+ {\rm sign} \ef)\geq 0.
\end{equation}
One sees that $\ef\in \mathcal V \Rightarrow (1\wedge \vert {\ef}\vert){\rm sign} \ef \in \mathcal V$ if and only if $f\in V_i\Longrightarrow(1\wedge |f|){\rm sign}f\in V_i$ for all $i=1,\ldots,m$. 
%
%
%
We have to prove the equivalence of the estimates in (ii) and \eqref{ouhcon}. Let first $\ef\in \mathcal C^\infty_i$. Then
$$(1\wedge |\ef|) {\rm sign} \ef=\begin{pmatrix}
f_1\\\vdots\\f_{j-1}\\(1\wedge |f_i|){\rm sign} {f_i}\\f_{j+1}\\ \vdots\\ f_m\end{pmatrix}\leftarrow i^{\rm th} \hbox{ row}$$
and 
$$(| \ef|-1)^+ {\rm sign} \ef=\begin{pmatrix}
	0\\\vdots\\ (|f_i|-1)^+{\rm sign} {f_i}\\ \vdots\\ 0\end{pmatrix}\leftarrow i^{\rm th} \hbox{ row}.$$ 
Accordingly, 
\begin{eqnarray*}
0&\leq& {\rm Re}{\ea}((1\wedge \vert \ef\vert){\rm sign} \ef,(\vert \ef\vert-1)^+ {\rm sign} \ef)\\
&=&\sum_{j\not=i} {\rm Re} a_{ij}(f_j,(|f_i|-1)^+{\rm sign} {f_i})+{\rm Re}a_{ii}((1\wedge|f_i|),(|f_i|-1)^+{\rm sign} {f_i})
\end{eqnarray*}
for all $\ef\in {\mathcal C}^\infty_i\cap {\mathcal V}$ and all $i=1,\ldots,m$. Due to the sesquilinearity of $a_{ij}$, this also implies
$$0\leq \sum_{j\not=i} {\rm Re}a_{ij}(\pm f_j,(|f_i|-1)^+{\rm sign} {f_i})+{\rm Re}a_{ii}((1\wedge|f_i|),(|f_i|-1)^+{\rm sign} {f_i})$$
for all $\ef\in {\mathcal C}^\infty_i\cap {\mathcal V}$, all $i=1,\ldots,m$, and all $\alpha\in{\mathbb C}$, $|\alpha|\leq 1$. This yields the claimed criterion. The converse implication can be proven analogously.
\end{proof}

\begin{rems}\label{xinftyrem}
(1) As we will see in Section~4, in many applications one has 
$$a_{ii}((1\wedge|f|){\rm sign}f,(|f|-1)^+{\rm sign} f)= 0\qquad\hbox{ for all }f\in V_i,\;i=1,\ldots,m.$$ 
In this case, it follows from the above theorem that a sufficient and necessary condition for $L^\infty$-contractivity of $(e^{t\ea})_{t\geq 0}$ is that for all $i=1,\ldots,m$
$$(1\wedge|f_i|){\rm sign}f_i \in V_i \qquad\hbox{and}\qquad a_{ij}(f_j,(|f_i|-1)^+{\rm sign} f_i)=0$$
for all $f_j\in B^\infty_{X_j}$, all $f_i\in V_i$, and all $j\not=i$.
This is a severe restriction to the possibility of extrapolating $(e^{t\ea})_{t\geq 0}$ to whole $L^p$-scale whenever our system ${\rm (ACP)}$ is actually coupled.

(2) The above result yields an alternative proof of~\cite[Lemma~6.1]{Mu07}. In fact, assume the spaces $V_i$ to have the following property: For each $f_j\in V_j$ one also has ${\rm sign}f_j\in V_j$. Then, after replacing $f_i$ by $f_i+{\rm sign}f_i$ in the condition in the above theorem, one sees that $(e^{t\ea})_{t\geq 0}$ is $L^\infty$-contractive if and only if
\begin{itemize}
\item $f\in V_i\Longrightarrow(1\wedge |f|){\rm sign}f\in V_i$ and
\item $\sum_{j\not= i} |a_{ij}(f_j,f_i)|\leq {\rm Re}a_{ii}({\rm sign}f_i,f_i)\hbox{ for all }\ef\in {\mathcal V}\cap {\mathcal C}^\infty_i.$
\end{itemize}
\end{rems}

As already mentioned, the main motivation for investigating $L^\infty$-contractivity is the extrapolation of $(e^{t\ea})_{t\geq 0}$ to $L^p$-spaces. 
We recall that if a semigroup $(T(t))_{t\geq 0}$ extrapolates to a consistent family of contractive $C_0$-semigroups on $L^p$, then it is called \emph{ultracontractive of dimension $d$} if there is a constant $c>0$ such that for all $p,q\in [1,\infty]$ and all $f\in L^p$ the estimate
$$\| T(t)f\|_{L^q}\leq ct^{-\frac{d}{2}|p^{-1}-q^{-1}|}\| f\|_{L^p}\qquad t\in[0,1],$$
holds, cf.~\cite[\S~7.3.2]{Ar04}.

\begin{thm}\label{ultra}
Let $H_i:=L^2(X_i)$. Assume $\ea$ to be accretive. Then the following assertions hold.
\begin{enumerate}[(1)]
\item The semigroup $(e^{t\ea})_{t\geq 0}$ extrapolates to a family of contractive $C_0$-semigroups on $L^p$, $p\in[1,\infty)$, which we denote again by $(e^{t\ea})_{\geq 0}$, if and only if for all $i=1,\ldots,m$ there holds
\begin{enumerate}[(i)]
\item $f\in V_i\Longrightarrow(1\wedge |f|){\rm sign}f\in V_i$,
\item $\sum_{j\not= i} |a_{ij}(f_j,(|f_i|-1)^+{\rm sign}f_i)|\leq {\rm Re}a_{ii}((1\wedge|f_i|){\rm sign}f_i,(|f_i|-1)^+{\rm sign} f_i)$ for all $\ef\in {\mathcal V}\cap {\mathcal C}^\infty_i$, and
\item $\sum_{j\not= i} |a_{ji}((|f_i|-1)^+{\rm sign}f_i,f_j)|\leq {\rm Re}a_{ii}((|f_i|-1)^+{\rm sign} f_i,(1\wedge|f_i|){\rm sign}f_i)$ for all $\ef\in {\mathcal V}\cap {\mathcal C}^\infty_i$.
\end{enumerate}
\item Let conditions (i)--(ii)--(iii) hold, assume that $V_i\cap L^1(X_i)$ is dense in $L^1(X_i)$, and let $d>2$ be a real number. Then $(e^{t\ea})_{t\geq 0}$ is \emph{ultracontractive of dimension $d$} if and only if $V_i$ is continuously embedded in $L^\frac{2d}{d-2}(X_i)$ for all $i=1,\ldots,m$.
\end{enumerate}
\end{thm}

\begin{proof}
(1) Let us first assume the semigroup $(e^{t\ea})_{t\geq 0}$ to extrapolate to a family of contractive $C_0$-semigroups on $L^p(X)$, $p\in[1,\infty)$, and hence in particular to be $L^\infty$-contractive. Moreover, since also the unit ball of $L^1(X)$ is left invariant, it follows by duality that the semigroup $(e^{t\ea^*})_{t\geq 0}$ is $L^\infty$-contractive. Here $\ea^*$ denotes the adjoint form of $\ea$, which by definition is given by
${\ea}^*(\ef,\eg)=\overline{{\ea}(\eg,\ef)}=\sum_{i,j=1}^m \overline{a_{ji}(g_i,f_j)}$, $\ef,\eg\in{\mathcal V}$. Since $\ea^*$ is accretive if and only if $\ea$ is accretive, we can apply Theorem~\ref{xinfty} to $(e^{t\ea})_{t\geq 0}$ and $(e^{t\ea^*})_{t\geq 0}$ and obtain conditions (i)--(ii)--(iii).

Conversely, since both $\ea$ and $\ea^*$ are accretive, it follows from (i)--(ii) and Theorem~\ref{xinfty} that $(e^{t\ea})_{t\geq 0}$ is $L^\infty$-contractive. Moreover, since also $\ea^*$ is accretive, it follows from (i)--(iii) and Theorem~\ref{xinfty} that $(e^{t\ea^*})_{t\geq 0}$ is $L^\infty$-contractive, too. Thus, by standard interpolation and duality methods one sees that $(e^{t\ea})_{t\geq 0}$ extrapolates to a family of contractive semigroups on $L^p(X)$, $p\in [1,\infty)$, that are strongly continuous for all $p>1$. Finally, by a result due to Voigt, contractivity implies strongly continuity of the extrapolated semigroup also in $L^1(X)$, cf.~\cite[\S~7.2.1]{Ar04}.

(2) The claim is a direct consequence of (1) and ~\cite[Thm.~7.3.2]{Ar04}.
\end{proof}

Observe that the extrapolated semigroups are positive in all $L^p$-spaces, $p\in [1,\infty)$, if and only if $(e^{t\ea})_{t\geq 0}$ is positive; they are analytic on all $L^p$-spaces, $p\in (1,\infty)$; finally, they are compact in all $L^p$-spaces, $p\in (1,\infty)$, if $V_i$ is compactly embedded in $L^2(X_i)$ for all $i=1,\ldots,n$. We refer the reader to~\cite[\S~7.3]{Ar04} for these and further properties of extrapolating semigroups.

\section{Applications}\label{applic}

\subsection{Ephaptical coupling of nerve fibres}\label{epha}
Motivated by the neurobiological theory of ephaptic coupling of myelinated fibres, cf. Remark~\ref{blbekbes}, we discuss a system
\begin{equation*}\label{problem2}
\left\{
\begin{array}{rcll}
\dot{u}_i(t,x)&=& \sum_{j=1}^m (c_{ij} u'_j(t,\cdot))'(x) &t\geq 0,\; x\in{\mathbb R},\;i=1,\ldots,m,\\
u_i(0,x)&=&u_{i0}(x), &x\in(0,1),\;i=1,\ldots,m,
\end{array}
\right.
\end{equation*}
of coupled diffusion equations on $m$ unbounded, parallel intervals. This case, which reflects the case of $m$ ephaptically interacting axons of infinite length, see e.g.~\cite{HK99},~\cite{BLBEK01}, and~\cite{BES01}, fits in the above framework if for $i,j=1,\ldots,m$ we let $H_i:=L^2({\mathbb R})$, $V_i:=H^1({\mathbb R})$, and
\begin{equation*}\label{ephaptform}
 a_{ij} (f,g):= \int_{-\infty}^{\infty} c_{ij}(x)f'(x) \overline{g'(x)} dx,\qquad f\in V_j,\;	g\in V_i.
\end{equation*}
It is known that~\eqref{problem2} is well-posed whenever the coefficients satisfy a uniform ellipticity condition, cf.~\cite{Am84}, and in fact the results of Section~\ref{sectionwellp} yield non-optimal criteria. However, assuming $\ea$ to be accretive we can perform an analysis of some qualitative properties of the system applying the theory developed in Section~\ref{sectionsymm}.
Let $\mathcal P$ be defined by
$$
\mathcal P\ef:=\left(\sum _{i=1}^m \frac{f_i}{m},\ldots,\sum _{i=1}^m \frac{f_i}{m}\right)^\top.
$$ 
Then $(e^{ta})_{t\geq 0}$ leaves invariant the closed subsets $\cya:=\left\{\ef\in {\mathcal H}:\|\ef-{\mathcal P}\ef \|\leq \alpha \right\}$, $\alpha\geq0$, if and only if there exist numbers $R_x \in \mathbb C,$ such that
\begin{equation}\label{primaephaptcond}
\sum _{j=1}^m c_{ij}(x)=R_x, \qquad\hbox{for a.e. } x\in{\mathbb R}, i=1,\ldots,m.
\end{equation}

Let now \eqref{primaephaptcond} hold. By Theorem~\ref{linainvar} the invariance of $\cya$ for some/all $\alpha\geq 0$ is equivalent to
\begin{equation}\label{boh}
\sum _{i,j=1}^m a_{ij}(g,h_i)=0, \qquad\hbox{for all } g,h_1,\ldots,h_m\in V \hbox{ s.t. }\sum_{i=1}^m h_i\equiv 0,
\end{equation}
i.e., to
$$ \sum_{i,j=i}^m \int_{-\infty}^{\infty} c_{ij}(x)g'(x) \overline{h_i'(x)} dx=0\quad \hbox{for all } g,h_1,\ldots,h_m\in V \hbox{ s.t. }\sum_{i=1}^m h_i\equiv 0.
$$
Since now $g'$ is indipendent of $i,j$ this is equivalent to
$$
\int_{-\infty}^{\infty} g'(x)\Big(\sum_{i,j=i}^m c_{ij}(x) \overline{h_i'(x)}\Big) dx=0\quad \hbox{for all }g,h_1,\ldots,h_m\in V \hbox{ s.t. }\sum_{i=1}^m h_i\equiv 0.
$$
Then, for a.e. $x\in\mathbb R$ there holds
$$
\sum_{i,j=i}^m c_{ij}(x) \overline{h_i'(x)}=\sum _{i=1}^m h'_{i}(x) \sum _{j=1}^m c_{ij}(x)= R_x \sum _{i=1}^m h'_{i}(x)=0\quad \hbox{for all }h_1,\ldots,h_m\in V \hbox{ s.t. }\sum_{i=1}^m h_i\equiv 0.
$$
 This shows that condition~\eqref{primaephaptcond} is sufficient. To see that it is also necessary, let~\eqref{boh} hold. 
Let $g,h \in H^1(\mathbb R)$ and consider the vector $\eh:=(h,-h,0,\ldots,0)^\top \in \mathcal V$. Since now $\sum _{i=1}^m h_i=0$ holds, we have
$$
\sum_{i,j=1}^m a_{ij}(g,h_i)= \int_{-\infty}^{\infty} g'(x)\left(\sum_{i,j=i}^m c_{ij}(x) \overline{h_i'(x)}\right) dx=0.
$$
Because of the arbitrarity of $g$ this yields that
$$
\sum_{i,j=i}^m c_{ij}(x) \overline{h_i'(x)}= h'(x)\sum _{j=1}^m c_{1j}(x)- h'(x)\sum _{j=1}^m c_{2j}(x)=0\qquad\hbox{for all }x\in\mathbb R.
$$
Iterating this procedure shows that~\eqref{boh} holds. Similarly, one shows that $(e^{ta})_{t\geq 0}$ leaves invariant the closed subsets $\bya:=\left\{\ef\in {\mathcal H}:\|{\mathcal P}\ef \|\leq \alpha \right\}$, $\alpha\geq 0$, if and only if there exist numbers $C_x \in \mathbb C,$ such that
\begin{equation*}\label{secondaephaptcond}
 \sum _{i=1}^m c_{ij}(x)=C_x, \qquad \hbox{a.e. for all }j=1,\ldots,m.
\end{equation*}

\begin{rem}\label{blbekbes}
Let us consider the case of two coupled axons. Several models based on~\eqref{problem2} have been proposed in the literature, assuming that $c_{11}=c_{22}=\alpha-\beta$ and $c_{21}=c_{12}=-\beta$, cf.~\cite[\S~4]{BES01}, or else $c_{11}=c_{22}=\alpha+\beta$, $c_{12}=c_{21}=-\beta$, cf.~\cite{BLBEK01}, or finally that $c_{11}=c_{12}=c_{21}=c_{22}$, cf.~\cite{HK99}, for some diffusion coefficient $\alpha$ and some coupling parameter $\beta>0$. Although these models are not equivalent,
%
%
in all of them the column and row sums agree, i.e., $c_{11} + c_{21}=c_{12} + c_{22}$ and $c_{11}+c_{12}=c_{21}+c_{22}$. Thus, condition~\eqref{invarmed} applies and the subsets $\{(f_1,f_2)\in L^2({\mathbb R})\times L^2({\mathbb R}):\|f_1-f_2\|_{L^2}\leq \alpha \}$ and $\{(f_1,f_2)\in L^2({\mathbb R})\times L^2({\mathbb R}):\|f_1+f_2\|_{L^2}\leq \alpha \}$ are left invariant for all $\alpha\geq 0$ under the action of $(e^{t\ea})_{t\geq 0}$.


\end{rem}

%
%

\subsection{A complete second order problem}\label{csop}

The strongly damped wave equation
\begin{equation}\label{problemwave}
\left\{
\begin{array}{rcll}
\ddot{u}(t,x)&=&\Delta(\alpha u+\dot{u})(t,x), &t\geq 0,\; x\in\Omega,\\
\frac{\partial u}{\partial \nu}(t,z)&=&\frac{\partial \dot{u}}{\partial \nu}(t,z)=0, &t\geq 0,\; z\in\partial \Omega,\\
u(0,x)&=&u_{0}(x), &x\in(0,1),\\
\dot{u}(0,x)&=&v_{0}(x), &x\in(0,1),
\end{array}
\right.
\end{equation}
on a bounded open domain of $\Omega\subset{\mathbb R}^n$, whose well-posedness has been proved in~\cite{Mu07b} for all $\alpha\in\mathbb C$, can also be treated with the methods presented in this paper.
%
The first equation has to be understood in the sense of distributions. 
We introduce a form $\ea:{\mathcal V}\times {\mathcal V}\to\mathbb C$, where 
 $$V_1=V_2=H_1=H^1(\Omega),\qquad H_2=L^2(\Omega),$$ 
 and
 \begin{eqnarray*}
 a_{11}(f,g)&:=&0,\qquad a_{12}(f,g):=-(f\mid g)_V=-\int_\Omega \nabla f\cdot \overline{\nabla g}dx,\\
 a_{21}(f,g)&:=&-\alpha \int_\Omega \nabla f\cdot \overline{\nabla g}dx,\qquad\hbox{and}\qquad 
 a_{22}(f,g):=\int_\Omega \nabla f\cdot \overline{\nabla g}dx.
 \end{eqnarray*}
Then $\ea$ is $\mathcal H$-elliptic and continuous due to Corollary~\ref{perturb0} and Proposition~\ref{continpure}, respectively, yielding well-posedness of~\eqref{problemwave}.
As an application of Proposition~\ref{productinv} we also mention that any closed subspace $\mathcal Y:=Y_1\times Y_2$ of the energy space $H^1(\Omega)\times L^2(\Omega)$  is invariant under the action of $(e^{t\ea})_{t\geq 0}$ if and only if
\begin{itemize}
\item $Y_2 \cap V \subset Y_1$,
\item $Y_2$ is invariant under the action of $(e^{ta_{22}})_{t\geq 0}$ for all $t \geq 0$, and
\item ${\rm Re}a_{21}(f,g)=0$ for all $f \in Y_1,$ $g \in Y^{2^\perp} \cap V$.
\end{itemize}
In this way one can e.g.~show that the solution $u$ to~\eqref{problemwave} has mean value $0$ as soon as the initial data $u_0,v_0$ have mean value $0$. Similarly, one can check that if $\Omega$ is a ball, then $u$ is a radial function provided that the initial data $u_0,v_0$ are radial. Such invariance properties of $(e^{t\ea})_{t\geq 0}$ directly follow from analogous ones of the Neumann heat semigroup $(e^{ta_{22}})_{t\geq 0}$. 


\subsection{A heat equation with dynamical boundary conditions}\label{dbc}

Let $\Omega$ be a bounded open domain of ${\mathbb R}^n$ with $C^\infty$ boundary $\partial\Omega$. Set
$$V_1:=H^1(\Omega),\quad H_1:=L^2(\Omega),\quad V_2:=H^1(\partial\Omega),\quad H_2:=L^2(\partial \Omega)$$
 and consider the initial-boundary value problem
\begin{equation}\label{cenn2}
\left\{
\begin{array}{rcll}
 \dot{u}(t,x)&=&\Delta u(t,x), &t\geq 0,\; x\in\Omega,\\
 \dot{w}(t,z)&=& u(t,z)+\Delta_{\partial \Omega} w(t,z), &t\geq 0,\;z\in\partial\Omega,\\
 w(t,z)&=&\frac{\partial u}{\partial \nu}(t,z), &t\geq 0,\; z\in\partial \Omega,\\
 u(0,x)&=&f(x), &x\in\Omega,\\
 w(0,z)&=&h(z), &z\in\partial\Omega.
\end{array}
\right.
\end{equation}
The results in this subsection should be compared with~\cite[\S~3]{CENN03} and~\cite[Ex.~5.6]{Mu06}, where well-posedness and exponential stability of~\eqref{cenn2} have also been investigated by different methods.
In~\eqref{cenn2} $\Delta_{\partial\Omega}$ denotes the Laplace--Beltrami operator, which is defined weakly as the operator associated with the form
$$a_{22}(f,g):=\int_{\partial\Omega} \nabla f\cdot \overline{\nabla g}d\sigma.$$
Moreover, define the forms
$$a_{11}(f,g):=\int_\Omega \nabla f\cdot \overline{\nabla g}dx,\qquad a_{12}(f,g):=-\int_{\partial \Omega} f\overline{g}_{|\partial \Omega} d\sigma,\qquad a_{21}(f,g):=-\int_{\partial \Omega} f_{|\partial \Omega}\overline{g} d\sigma.$$
A direct integration by parts shows that the operator $\mathcal A$ associated with $\ea$ is the same one that governs the above problem, i.e.,
$${\mathcal A}:=\begin{pmatrix}
\Delta & 0\\
\cdot_{|\partial \Omega} & \Delta_{\partial \Omega}\end{pmatrix},\qquad D({\mathcal A})=\left\{\begin{pmatrix}u\\ w\end{pmatrix}\in H^2(\Omega)\times H^2(\partial \Omega): \frac{\partial u}{\partial \nu}=w\right\}.$$
We show that $\mathcal A$ generates a semigroup that is analytic of angle $\frac{\pi}{2}$ and positive, but which does not leave the unit ball of $L^\infty(\Omega)\times L^\infty(\partial \Omega)$ invariant.

First, observe that $a_{11}$ (resp. $a_{22}$) are continuous and $H_1$- (resp. $H_2$-) elliptic. Moreover, due to boundedness from $H^1(\Omega)$ to $L^2(\partial\Omega)$ of the trace operator, forms $a_{12}$ and $a_{21}$ are bounded on $H_2\times V_1$ and on $H_1\times V_2$, respectively. Accordingly, $\ea$ is continuous and by Proposition~\ref{perturb0} also $H_1\times H_2$-elliptic.
In order to apply Proposition~\ref{cosine}, observe that ${\rm Im}a_{ii}(f,f)=0$ for all $f\in V_i$, $i=1,2$. Moreover
\begin{eqnarray*}
|{\rm Im}\left(a_{12}(f,g)+a_{21}(g,f)\right)|&=&|{\rm Im}(\int_{\partial \Omega} f\overline{g}_{|\partial \Omega} d\sigma+\int_{\partial \Omega} g_{|\partial \Omega}\overline{f} d\sigma)|\\
&=&|{\rm Im}(\int_{\partial \Omega} f\overline{g}_{|\partial \Omega} d\sigma+\overline{\int_{\partial \Omega} f\overline{g}_{|\partial \Omega} d\sigma})|=0.
\end{eqnarray*}

By Theorem~\ref{contractive}, the semigroup is real. To see that it is positive, observe that $a_{11}$ is associated with the Laplace operator with Neumann boundary conditions and $a_{22}$ with the Laplace--Beltrami operator on $\partial\Omega$. Therefore they generate positive semigroups and the first condition of Theorem~\ref{contractive}.(2) is satisfied. The second condition is also clear since $f_{|\partial \Omega}$ is positive whenever $f$ is positive. By Corollary~\ref{domindiag}, $(e^{t\ea})_{t\geq 0}$ dominates the semigroup $(e^{t\ea_0})_{t\geq0}$, where $\ea_0:=a_{11}+a_{22}$, which governs the uncoupled system of two diffusion equations on $\Omega$ (with homogeneous Neumann boundary conditions) and $\partial\Omega$. 

It also follows from Theorem~\ref{xinfty} and Remark~\ref{xinftyrem}.(1) that $(e^{t\ea})_{t\geq 0}$ is not $L^\infty(\Omega)\times L^\infty(\partial\Omega)$-contractive, since for non-constant $f\in H^1(\partial\Omega)$ such that $|f|\leq 1$ and for $g\in H^1(\Omega)$ with $g_{|\partial\Omega}=1+f$ one has $a_{12}(f,g)=-\int_{\partial \Omega} |\nabla f|^2 d\sigma<0$, which contradicts condition (ii) in Theorem~\ref{xinfty}. However, Theorem~\ref{ultra} can be used in order to show $L^p$-well-posedness for~\eqref{cenn2}.

We first prove a generation result in all $L^p$-spaces for $p\geq 2$. Write $\mathcal A$ as
$${\mathcal A}:=\tilde{\mathcal A}+{\mathcal B}:=\begin{pmatrix}
\Delta-C^* & 0\\
\cdot_{|\partial \Omega} & \Delta_{\partial \Omega}-Id\end{pmatrix}+\begin{pmatrix}
C^* & 0\\
0 &Id \end{pmatrix}.$$
Here, $C^*$ is the adjoint of the linear operator from $H^1(\Omega)$ to $L^2(\Omega)$ defined by
$$(Cf)(x):=\nabla f(x) \cdot \overline{\nabla D_N1(x)},\qquad f\in H^1(\Omega),\; x\in\Omega,$$
where $D_N1$ denotes the unique (modulo constants) solution $u$ of 
$$\left\{
\begin{array}{rcll}
\Delta u(x)&=& 0,\qquad &x\in\Omega,\\
\frac{\partial u}{\partial \nu}(z)&=&1,&z\in\partial \Omega.
\end{array}
\right.$$
The operator $\tilde{\mathcal A}$ is associated with the matrix form $\tilde{\ea}$ whose entries are given by
\begin{eqnarray*}
\tilde{a}_{11}(f,g)&:=&\int_\Omega \nabla f\cdot \overline{\nabla g}dx+\int_\Omega f \overline{(\nabla g\cdot \nabla D_N 1)} dx,\\ \tilde{a}_{22}(f,g)&:=&\int_{\partial\Omega} \nabla f \cdot \overline{\nabla g}d\sigma +\int_{\partial\Omega}f\overline{g}d\sigma,\quad\hbox{and}\\
\tilde{a}_{12}&:=&a_{12}\qquad\hbox{as well as}\qquad \tilde{a}_{21}\quad:=\quad a_{21}.
\end{eqnarray*}
One sees that the perturbation $\tilde{a}_{11}-a_{11}$ is bounded on $H_1\times V_1$, thus by Lemma~\ref{perturbd} $\tilde{\ea}$ is associated with a semigroup $(e^{t\tilde{\ea}})_{t\geq 0}$ on $H_1\times H_2$. For all $g\in H^1(\partial\Omega)$ such that $|g|\leq 1$ and all $f\in H^1(\Omega)$ we have
\begin{eqnarray*}
|\tilde{a}_{12}(g,(|f|-1)^+ {\rm sign}{f})|
&\leq &\int_{\partial\Omega} |g| (|f|-1)^+ d\sigma\leq \int_{\partial\Omega} (|f|-1)^+ d\sigma\\
&=&\int_{\partial\Omega} \frac{\partial D_N 1}{\partial \nu}(|f|-1)^+ d\sigma\\
&= & \int_{\Omega} \nabla(|f|-1)^+\cdot \nabla D_N 1 dx\\
&=&\int_\Omega (1\wedge |f|) \left(\nabla(|f|-1)^+ \cdot \nabla D_N 1\right) {\mathbb 1}_{\{|f|\geq 1\}}dx\\
&&\quad +\int_\Omega (1\wedge |f|) \left(\nabla(|f|-1)^+ \cdot \nabla D_N 1\right) {\mathbb 1}_{\{|f|\leq 1\}}dx\\
&=&{\rm Re}\tilde{a}_{11} ((1\wedge |f|){\rm sign}f,(|f|-1)^+{\rm sign}f).
\end{eqnarray*}
since $\nabla(|f|-1)^+=0$ a.e. on $\{x\in \Omega : |f(x)|\leq 1\}$. Likewise, for all $f\in H^1(\Omega)$ such that $|f|\leq 1$ (so that in particular $|f_{|\partial\Omega}|\leq 1$) and all $g\in H^1(\partial\Omega)$
\begin{eqnarray*}
|\tilde{a}_{21}(f,(|g|-1)^+ {\rm sign}{g})|
&\leq &\int_{\partial\Omega} |f| (|g|-1)^+ d\sigma\leq \int_{\partial\Omega} (|g|-1)^+ d\sigma\\
&=& \int_{\partial\Omega} (1\wedge |g|) (|g|-1)^+ {\mathbb 1}_{\{|g|\geq 1\}}dx\\
&&\quad +\int_{\partial\Omega}(1\wedge |g|) (|g|-1)^+ {\mathbb 1}_{\{|g|\leq 1\}} dx\\
&=&{\rm Re}\tilde{a}_{22} ((1\wedge |g|){\rm sign}g,(|g|-1)^+{\rm sign}g).
\end{eqnarray*}
Thus, Theorem~\ref{xinfty} applies and we conclude that $(e^{t\tilde{\ea}})_{t\geq 0}$ extrapolates to a consistent family of semigroups on $L^p(\Omega)\times L^p(\partial\Omega)$, $p\geq  2$, the generator of the semigroup in $L^p(\Omega)\times L^p(\partial \Omega)$ being the part of $\tilde{\mathcal A}$ in $L^p(\Omega)\times L^p(\partial \Omega)$. Since now (the part of) $\mathcal B$ is compact from $W^{2,p}(\Omega)\times W^{2,p}(\partial \Omega)$ to $L^p(\Omega)\times L^p(\partial\Omega)$ for all $p=[1,\infty)$, by the perturbation thorem of Desch--Schappacher (see e.g.~\cite[Thm.~3.7.25]{ABHN01}) we conclude that (the part of) ${\mathcal A}=\tilde{\mathcal A}+{\mathcal B}$ generates a semigroup on $L^p(\Omega)\times L^p(\partial\Omega)$, $p\geq 2$.

Introducing a different operator $\tilde{\mathcal A}$ (more precisely, replacing $C^*$ by $C$), and hence a different perturbed form $\tilde{\ea}$, it is also possible to prove in a similar manner that the semigroup associated with the adjoint $\tilde{\ea}^*$ is $L^\infty$-contractive. By duality we conclude as above that $\mathcal A$ generates a semigroups on $L^p(\Omega)\times L^p(\partial\Omega)$ also for $p\in [1,2]$, hence for the whole scale of $L^p$-spaces. 

Observe, however, that none of these semigroups is ultracontractive. For example, in the first considered case, $|\tilde{a}_{12}((|g|-1)^+ {\rm sign}{g},f)|\leq {\rm Re}\tilde{a}_{22} ((|g|-1)^+{\rm sign}g,(1\wedge |g|){\rm sign}g)$ does not hold for all $(f,g)\in H^1(\Omega)\times H^1(\partial\Omega)$ such that $|f|\leq 1$, hence condition (iii) in Theorem~\ref{ultra}.(1) is not satisfied. Thus we {cannot} deduce ultracontractivity of $(e^{t\tilde{\ea}})_{t\geq 0}$ 
from  Theorem~\ref{ultra}.(2) and the Sobolev embeddings $H^1(\Omega)\hookrightarrow L^\frac{2n}{n-2}(\Omega)$, $H^1(\partial\Omega)\hookrightarrow L^\frac{2n-2}{n-3}(\partial\Omega)$.

\section{Appendix: Hilbert space projections}

In the following we denote by $\mathcal P$ the orthogonal projection onto a closed subspace $Y$ of a Hilbert space $\mathcal H$, and by $\cya$ the closed convex subset of $\mathcal H$ defined as the strip around $Y$ of thickness $2\alpha$, i.e.,
$$
\cya:=\left\{f \in {\mathcal H} : \|f -{\mathcal P}f \|\leq \alpha \right\}.
$$
A subset ${\mathcal S}\subset \mathcal H$ is said to be \emph{invariant under a semigroup $({\mathcal T}(t))_{t\geq 0}$}, if ${\mathcal T}(t){\mathcal S}\subset {\mathcal S}$ for all $t\geq 0$.

\begin{prop}\label{ficata}
Let $({\mathcal T}(t))_{t\geq 0}$ be a $C_0$-semigroup on $\mathcal H$. Consider the following assertions.
\begin{itemize}
\item[\emph{(a)}] $\cya$ is invariant under $({\mathcal T}(t))_{t\geq 0}$ for all $\alpha>0$.
\item[\emph{(b)}] $\cyb$ is invariant under $({\mathcal T}(t))_{t\geq 0}$ for some $\beta>0$.
\item[\emph{(c)}] $Y$ is invariant under $({\mathcal T}(t))_{t\geq 0}$.
\end{itemize}
Then (a)$ \Longleftrightarrow $(b)$ \Longrightarrow $(c). If $({\mathcal T}(t))_{t\geq 0}$ is contractive, then also (c) $ \Longrightarrow$ (a) holds.
\end{prop}

\begin{proof}
$(a) \Longrightarrow (b)$ is trivial. In order to prove the converse implication, observe that ${\mathcal C}_{Y,\beta}=\frac{\beta}{\alpha}\cya$ for all $\beta>0$, since $\mathcal P$ is linear. The claim follows from linearity of $({\mathcal T}(t))_{t\geq 0}$.


In order to prove $(a) \Longrightarrow (c)$, let $f\in Y$. Thus, $\| f -{\mathcal P}f \|\leq \alpha$ for all $\alpha>0$. Since $({\mathcal T}(t))_{t\geq 0}$ leaves invariant $\cya$, $\| {\mathcal T}(t)f-{\mathcal P}{\mathcal T}(t)f\|\leq \alpha$ for all $\alpha>0$, i.e., ${\mathcal T}(t)f={\mathcal P}{\mathcal T}(t)f$ and ${\mathcal T}(t)f\in Y$ for all $t\geq 0$.


Finally, let $({\mathcal T}(t))_{t\geq 0}$ be contractive. To this aim, let $f \in\cyb$ and observe that there exists $f_0 \in Y$ such that $\|f-f_0\|\leq \beta$. Furthermore, due to the contractivity of $({\mathcal T}(t))_{t\geq 0}$ one has $\|{\mathcal T}(t)f-{\mathcal T}(t)f_0\|\leq \beta$. Since $T(t)$ leaves $Y$ invariant, one has $T(t)f_0 \in Y$. Since ${\mathcal P}{\mathcal T}(t)y$ is the element of $Y$ with minimal distance from ${\mathcal T}(t)y$, we conclude that $\| {\mathcal T}(t)f-{\mathcal P}{\mathcal T}(t)f\|\leq \|{\mathcal T}(t)f-{\mathcal T}(t)f_0\|\leq \beta$, i.e., ${\mathcal T}(t)y\in \cyb$.
\end{proof}

In the special case of semigroups coming from a sesquilinear form we obtain the following.

\begin{cor}\label{lemmino}
If the semigroup $({\mathcal T}(t))_{t\geq 0}$ is associated with a densely defined, ${\mathcal H}$-elliptic, continuous form $a:V\times V\to{\mathbb C}$, then the following assertions are equivalent.
\begin{itemize}
\item[\emph{(c)}] $Y$ is invariant under $({\mathcal T}(t))_{t\geq 0}$.
\item[\emph{(d)}] For all $f\in V$, $g\in{\rm ker}(I-{\mathcal P})\cap V$, and $h\in{\rm ker}{\mathcal P}\cap V$ there holds ${\mathcal P}f\in V$ and $a(g, h)= 0$. 
\end{itemize}
If in particular $a$ is accretive, then the assertions (a)--(d) above are all equivalent.
\end{cor}

\begin{proof}
Under the above assumptions we can apply~\cite[Thm.~2.2]{Ou04} and directly obtain that (a)--(c) in Proposition~\ref{ficata} are equivalent to the condition that {for all } $f\in V$ there holds ${\mathcal P}f\in V$   and ${\rm Re}a({\mathcal P}f,f-{\mathcal P}f)\geq 0$. Observe that accretivity of $a$, which is an assumption of~\cite[Thm.~2.2]{Ou04}, is not needed while studying invariance of {subspaces}. Taking into account the decomposition ${\mathcal H}={\rm ker}(I-{\mathcal P})\oplus {\rm ker}{\mathcal P}$ and the sesquilinearity of $a$
we obtain that the above condition is equivalent to the claimed criterion.
\end{proof}

Similarly, since $I-{\mathcal P}$ is the projection of ${\mathcal H}$ onto ${\rm ker}\;{\mathcal P}$, then the following also holds. 

\begin{cor}\label{lemmino2}
Let $a:V\times V\to{\mathbb C}$ be a densely defined, accretive, ${\mathcal H}$-elliptic, continuous form. The following assertions are equivalent, where we use the notation $\bya:=\left\{f \in {\mathcal H} : \|{\mathcal P}f \|\leq \alpha \right\}$.
\begin{itemize}
\item[\emph{(a)}] $\byb$ is invariant under $({\mathcal T}(t))_{t\geq 0}$ for some $\beta>0$.
\item[\emph{(b)}] $\bya$ is invariant under $({\mathcal T}(t))_{t\geq 0}$ for all $\alpha>0$.
\item[\emph{(c)}] $Y$ is invariant under $({\mathcal T}(t))_{t\geq 0}$.
\item[\emph{(d)}] For all $f \in \mathcal V$, $g\in{\rm ker}(I-{\mathcal P})\cap V$, and $h\in{\rm ker}{\mathcal P}\cap V$ there holds ${\mathcal P}f \in V$ and $a(h,g)= 0$.
\end{itemize}
\end{cor}

\begin{rem}
%
Corollary~\ref{lemmino} directly follows from~\cite[Thm.~2.2]{Ou04}. It yields new proofs of known facts. E.g., let $H=L^2(\Omega)$, $\Omega$ an open ball, and $Y$ be the space of radial functions over $\Omega$, i.e. of functions $f$ such that $f(x)=f(y)$ if $|x|=|y|$. Then switching to polar coordinates and applying Fubini's theorem yields that $Y$ is invariant under the semigroup generated by the Laplacian with several ``radial" boundary conditions -- including Dirichlet and Neumann ones. Provided that the boundary coefficient is constant, one can show in this way that radial initial value give to radial solutions also under Robin, Wentzell--Robin, and those dynamical boundary conditions considered in \S~4.3.\end{rem}.

\end{document}